\documentclass[aos,seceqn,citesort,dvips]{arximspdf}
\usepackage{newsym2e}
\usepackage{dcolumn}
\usepackage{graphicx}

\doi{10.1214/08-AOS663}
\volume{37}
\issue{6A}
\pubyear{2009}
\firstpage{3697}
\lastpage{3714}

\makeatletter

\newcolumntype{d}[1]{D{.}{.}{#1}}

\newcommand{\mesf}{\mathsf}

\setattribute{copyright}{owner}{\textup{In the Public Domain}}

\makeatother

\begin{document}
\begin{frontmatter}

\title{Estimating the Gumbel scale parameter for local alignment of random
sequences by importance sampling with stopping times}
\runtitle{Gumbel scale parameter estimation for local alignment}

\begin{aug}
\author[A]{\fnms{Yonil} \snm{Park}\ead[label=e1]{park@ncbi.nlm.nih.gov}},
\author[A]{\fnms{Sergey} \snm{Sheetlin}\ead[label=e2]{sheetlin@ncbi.nlm.nih.gov}} and
\author[A]{\fnms{John L.} \snm{Spouge}\corref{}\ead[label=e3]{spouge@ncbi.nlm.nih.gov}}
\runauthor{Y. Park, S. Sheetlin and J. L. Spouge}
\affiliation{National Library of Medicine}
\address[A]{Y. Park\\
S. Sheetlin\\
J. L. Spouge\\
National Center\\
\quad for Biotechnology Information\\
National Library of Medicine\\
National Institutes\\
\quad of Health\\
8600 Rockville Pike\\
Bethesda, Maryland 20894\\
USA\\
\printead{e1}\\
\phantom{E-mail: }\printead*{e2}\\
\phantom{E-mail: }\printead*{e3}} 
\end{aug}

\received{\smonth{10} \syear{2007}}
\revised{\smonth{6} \syear{2008}}

\begin{abstract}
The gapped local alignment score of two random sequences follows
a~Gumbel distribution. If computers could estimate the parameters of the
Gumbel distribution within one second, the use of arbitrary alignment
scoring schemes could increase the sensitivity of searching biological
sequence databases over the web. Accordingly, this article gives a novel
equation for the scale parameter of the relevant Gumbel distribution. We
speculate that the equation is exact, although present numerical
evidence is limited. The equation involves ascending ladder variates in
the global alignment of random sequences. In global alignment
simulations, the ladder variates yield stopping times specifying random
sequence lengths. Because of the random lengths, and because our trial
distribution for importance sampling occurs on a different sample space
from our target distribution, our study led to a mapping theorem, which
led naturally in turn to an efficient dynamic programming algorithm for
the importance sampling weights. Numerical studies using several popular
alignment scoring schemes then examined the efficiency and accuracy of
the resulting simulations.
\end{abstract}


\end{frontmatter}

\section{Introduction}\label{sec1}

Sequence alignment is an indispensable tool in modern molecular biology.
As an example, BLAST \cite{1,2,3} (the Basic Local Alignment Search Tool,
\url{http://www.ncbi.nlm.nih.gov/BLAST/}), a popular sequence alignment
program, receives about $2.89$ submissions per second over the Internet.
Currently, BLAST users can choose among only 5 standard alignment
scoring systems, because BLAST $p$-values must be pre-computed with
simulations that take about 2 days for the required $p$-value
accuracies. Moreover, adjustments for unusual amino acid compositions
are essential in protein database searches \cite{4}, and in that application,
computational speed demands that the corresponding $p$-values be
calculated with crude, relatively inaccurate approximations \cite{2}.
Accordingly, for more than a decade, much research has been directed at
estimating BLAST $p$-values in real time (i.e., in less than 1 sec)
\cite{5,6,7,8}, so that BLAST might use arbitrary alignment scoring systems.

Several studies have used importance sampling to estimate the BLAST
$p$-value \cite{6,8,9}. To describe importance sampling briefly, let
$\mathbb{E}$ denote the expectation for some ``target
distribution'' $\mathbb{P}$, let $\mathbb{Q}$ be any distribution, and
consider the equation
%
\begin{equation}\label{eq01}
\mathbb{E}X: = \int X( \omega)\,d\mathbb{P}( \omega) = \int X(
\omega)\frac{d\mathbb{P}( \omega )}{d\mathbb{Q}( \omega )}\,d\mathbb{Q}(
\omega).
\end{equation}

A computer can draw samples $\omega_{i}$ $(i = 1,\ldots,r)$ from the ``trial
distribution'' $\mathbb{Q}$ to estimate the expectation: $\mathbb{E}X
\approx r^{ - 1}\sum_{i = 1}^{r} X( \omega_{i} )[ d\mathbb{P}(
\omega_{i} )/d\mathbb{Q}( \omega_{i} ) ]$. The name ``importance
sampling'' derives from the fact that the subsets of the sample space
where~$X$ is large dominate contributions to $\mathbb{E}X$. By focusing
sampling on the ``important'' subsets, judicious choice of the trial
distribution $\mathbb{Q}$ can reduce the effort required to estimate
$\mathbb{E}X$. In importance sampling, the likelihood ratio
$d\mathbb{P}( \omega)/d\mathbb{Q}( \omega)$ is often called the
``importance sampling weight'' (or simply, the ``weight'') of the
sample $\omega$.

A Monte Carlo technique called ``sequential importance sampling'' can
substantially increase the statistical efficiency of importance sampling
by generating samples from $\mathbb{Q}$ incrementally and exploiting the
information gained during the increments to guide further increments.
Although sequences might seem an especially natural domain for
sequential sampling, most simulation studies for BLAST $p$-values have
used sequences of fixed length. In contrast, our study involves
sequences of random length.

Here, as in several other importance sampling studies \cite{6,8,9,10}, hidden
Markov models generate a trial distribution $\mathbb{Q}$ of random
\textit{alignments} between two sequences, where the
\textit{sequences} have a target distribution $\mathbb{P}$. The other
studies gloss over the fact that their trial and target distributions
occur on different sample spaces, such as alignments and sequences. The
other studies used sequences of fixed lengths, however, where a
relatively simple formula for the weight $d\mathbb{P}/d\mathbb{Q}$
pertains. For the sequences of random length in this paper, however, the
stopping rules for sequential sampling complicate formulas
for $d\mathbb{P}/d\mathbb{Q}$. Accordingly, the \hyperref[app]{Appendix} gives a general
mapping theorem giving formulas for the weights
$d\mathbb{P}/d\mathbb{Q}$ when each sample from $\mathbb{P}$ corresponds
to many different samples from $\mathbb{Q}$. (In the present article,
e.g., each pair of random sequences corresponds to many possible random
alignments.) In addition to the mapping theorem, we also develop several
other techniques specifically tailored to speeding the estimation of the
BLAST $p$-value.

The organization of this article follows. Section \ref{sec2} on background and
notation is divided into 4 subsections containing: (1) a friendly
introduction to sequence alignment and its notation; (2) a brief
self-contained description of the algorithm for calculating global
alignment scores; (3) a technical summary of previous research on
estimating the BLAST $p$-value introducing our importance sampling
methods; and (4) a heuristic model for random sequence alignment using
Markov additive processes. Section \ref{sec3} on Methods is also divided into 4
subsections containing: (1) a novel formula for the relevant Gumbel
scale parameter $\lambda$; (2) a Markov chain model for simulating
sequence alignments (borrowed directly from a previous study \cite{10}, but
used here with a stopping time); (3) a dynamic programming algorithm for
calculating the importance sampling weights in the presence of a
stopping time; and (4) formulas for the simulation errors. Section \ref{sec4}
then gives numerical results for the estimation of $\lambda$ under 5
popular alignment scoring schemes. Finally, Section \ref{sec5} is our Discussion.

\section{Background and notation}\label{sec2}

\subsection{Sequence alignment and its notation}\label{sec21}

Let $\mathbf{A} = A_{1}A_{2}\cdots$ and $\mathbf{B} = B_{1}B_{2}\cdots$ be two
semi-infinite sequences drawn from a finite alphabet $\mathfrak{L}$, for example,
$\{ \mathrm{A},\mathrm{C},\mathrm{D},\mathrm{E},\mathrm{F},\mathrm{G},\mathrm{H},
\mathrm{I},\mathrm{K},\mathrm{L},\mathrm{M},\mathrm{N},\mathrm{P},\mathrm{Q},
\mathrm{R},\mathrm{S},\mathrm{T},\mathrm{V},\mathrm{W},\mathrm{Y} \}$ (the amino acid
alphabet) or $\{ \mathrm{A},\mathrm{C},\mathrm{G},\mathrm{T} \}$ (the nucleotide alphabet). Let
$s\dvtx \mathfrak{L} \times \mathfrak{L} \mapsto\mathbb{R}$ denote a ``scoring matrix.'' In
database applications, $s( a,b )$ quantifies the similarity between $a$
and~$b$, for example, the so-called ``PAM'' (point accepted mutation) and
``\mbox{BLOSUM}'' (block sum) scoring matrices can quantify evolutionary
similarity between two amino acids \cite{11,12}.

\begin{figure}

\includegraphics{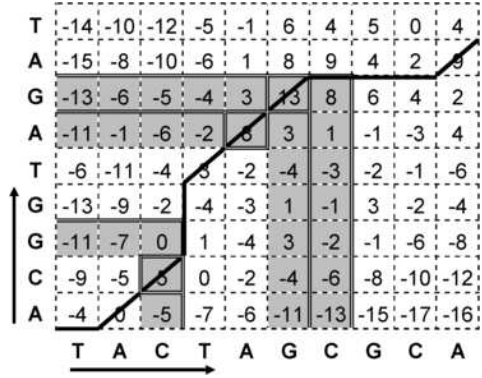}

\caption{Gapped global alignment scores and the corresponding directed
paths for two subsequences $\mathbf{A}[ 1,10 ] = \mathsf{TACTAGCGCA}$
and $\mathbf{B}[ 1,9 ] = \mathsf{ACGGTAGAT}$, drawn from the nucleotide
alphabet $\{ \mathsf{A},\mathsf{C},\mathsf{G},\mathsf{T} \}$. Figure \protect\ref{fig1} uses a nucleotide scoring matrix,
where $s( a,b ) = 5$ if $a = b$ and $- 4$ otherwise, and the affine gap
penalty $w_{g} = 3 + 2g$. The vertex $( i,j )$ is in the northeast corner
of the cell $( i,j )$, with the origin $( 0,0 )$ at the southwest corner
of Figure \protect\ref{fig1}. The cell $( i,j )$ displays the global score~$S_{i,j}$,
calculated from (\protect\ref{eq12}). The optimal global path ending at the
point $( 10,8 )$, for example, consists of $12$ edges, in order: $1$ east
of length $1, 2$ northeast, $1$ north of length $2, 3$ northeast, $1$ east
of length $3$, and $1$ northeast. The optimal global score $S_{10,8} = -
5 + 5 + 5 - 7 + 5 + 5 + 5 - 9 + 5 = 9$ is the sum of the corresponding
edges and represents the path of greatest weight starting at $( 0,0 )$
and ending at $( 10,8 )$. The corresponding optimal global alignment of
the subsequences $\mathbf{A}[ 1,10 ]$ and $\mathbf{B}[ 1,9 ]$ is\break
\mbox{\hspace*{148pt}$\mathsf{TAC} \mbox{--}\, \mbox{--} \mathsf{TAGCGCA}$} \break
\mbox{\hspace*{146pt}$- \mathsf{ACGGTAG} \mbox{--}\, \mbox{--}\, \mbox{--} \mathsf{A}$}.\break
The edge maxima are $M_{1} = - 4, M_{2} = 0, M_{3} = 5,M_{4} = 1 , M_{5}
= 3, M_{6} = 8, M_{7} = 13, M_{8} = 9$, $M_{9} = 6$. The shading and the
double lines indicate squares where a vertex (surrounded by double
lines) generated an SALE $\beta( k )$. The SALE scores are $M_{\beta ( 1
)} =\break M_{3} = 5, M_{\beta ( 2 )} = M_{6} = 8,M_{\beta ( 3 )} = M_{7} =
13;$ and the global maximum $M$ for $\mathbf{A}$ and $\mathbf{B}$ is no
less than $13$, the largest global score shown.}
\label{fig1}
\end{figure}

The alignment graph $\Gamma_{\mathbf{A},\mathbf{B}}$ of the
sequence-pair $( \mathbf{A},\mathbf{B} )$ is a directed, weighted
lattice graph in two dimensions, as follows. The vertices $v$ of
$\Gamma_{\mathbf{A},\mathbf{B}}$ are nonnegative integer points $( i,j
)$. (Below, ``$: =$'' denotes a definition, e.g., the natural numbers
are $\mathbb{N}: = \{ 1,2,3,\ldots\}$. Throughout the article, $i, j, k, m,
n$ and $g$ are integers.) Three sets of directed edges $e$ come out of
each vertex $v = ( i,j )$: northward, northeastward and eastward (see
Figure \ref{fig1}). One northeastward edge goes into $v = ( i + 1,j + 1 )$ with
weight $s[ e ] = s( A_{i + 1},B_{j + 1} )$. For each $g > 0$, one eastward
edge goes into $v = ( i + g,j )$ and one northward edge goes into $v = (
i,j + g )$; both are assigned the same weight $s[ e ] = - w_{g} < 0$. The
deterministic function $w\dvtx\mathbb{N} \mapsto( 0,\infty]$ is called the
``gap penalty.'' (The value $w_{g} = \infty$ is explicitly permitted.)
This article focuses on affine gap penalties $w_{g} = \Delta_{0} +
\Delta_{1}g$ $(\Delta_{0},\Delta_{1} \ge 0)$, which are typical in BLAST
sequence alignments. Together, the scoring matrix $s( a,b )$ and the gap
penalty $w_{g}$ constitute the ``alignment parameters.''

A (directed) path $\pi = ( v_{0},e_{1},v_{1},e_{2},\ldots,e_{k},v_{k} )$ in
$\Gamma_{\mathbf{A},\mathbf{B}}$ is a finite alternating sequence of
vertices and edges that starts and ends with a vertex. For each $i =
1,2,\ldots,k$, the directed edge $e_{i}$ comes out of vertex $v_{i - 1}$
and goes into vertex $v_{i}$. We say that the path $\pi$ starts at
$v_{0}$ and ends at $v_{k}$.

Denote finite subsequences of the sequence $\mathbf{A}$ by $\mathbf{A}[
i,m ] = A_{i}A_{i + 1}\cdots A_{m}$. Every gapped alignment of the
subsequences $\mathbf{A}[ i,m ]$ and $\mathbf{B}[ j,n ]$ corresponds to
exactly one path that starts at $v_{0} = ( i - 1,j - 1 )$ and ends at
$v_{k} = ( m,n )$ (see Figure \ref{fig1}). The alignment's score is the ``path
weight'' $S_{\pi} : = \sum_{i = 1}^{k} s[ e_{i} ]$.

Define the ``global score'' $S_{i,j}: = \max_{\pi} S_{\pi}$, where the
maximum is taken over all paths $\pi$ starting at $v_{0} = ( 0,0 )$ and
ending at $v_{k} = ( i,j )$. The paths $\pi$ starting at $v_{0}$, ending
at $v_{k}$, and having weight $S_{\pi} = S_{i,j}$ are ``optimal global
paths'' and correspond to ``optimal global alignments'' between
$\mathbf{A}[ 1,i ]$ and $\mathbf{B}[ 1,j ]$. Define the ``edge
maximum'' $M_{n}: = \max\{ \max_{0 \le i \le n}S_{i,n},\max_{0 \le j \le
n}S_{n,j} \}$, and the ``global maximum'' $M: = \sup_{n \ge 0}M_{n}$.
(The single subscript in $M_{n}$ indicates that the variate corresponds
to a square $[ 0,n ] \times[ 0,n ]$, rather than a general rectangle $[
0,m ] \times[ 0,n ]$.) Define the ``strict ascending ladder epochs''
(SALEs) in the sequence $( M_{n} )$: let $\beta( 0 ): = 0$ and $\beta( k +
1 ): = \min\{ n > \beta( k )\dvtx M_{n} > M_{\beta ( k )} \}$,
where $\min\varnothing : = \infty$. We call $M_{\beta ( k )}$ the ``$k$th
SALE score.''

Define also the ``local score'' $\tilde{S}_{i,j}: = \max_{\pi} S_{\pi}$,
where the maximum is taken over all paths $\pi$ ending at $v_{k} = ( i,j
)$, regardless of their starting point. Define the ``local
maximum'' $\tilde{M}_{m,n}: = \max_{0 \le i \le m,0 \le j \le
n}\tilde{S}_{i,j}$. The paths $\pi$ ending at $v_{k} = ( i,j )$ with
local score $S_{\pi} = \tilde{S}_{i,j} = \tilde{M}_{m,n}$ are ``optimal
local paths'' corresponding to the ``optimal local alignments'' between
subsequences of $\mathbf{A}[ 1,m ]$ and $\mathbf{B}[ 1,n ]$.

Now, the following ``independent letters'' model introduces randomness.
Choose each letter in the sequence $\mathbf{A}$ and $\mathbf{B}$
randomly and independently from the alphabet $\mathfrak{L}$ according to
fixed probability distributions $\{ p_{a}\dvtx a \in\mathfrak{L} \}$ and $\{
p'_{b}\dvtx b \in\mathfrak{L} \}$. (Although this article permits the
distributions $\{ p_{a} \}$ and $\{ p'_{b} \}$ to be different, in
applications they are usually the same.) Throughout the paper, the
probability and expectation for the independent letters model are
denoted by $\mathbb{P}$ and $\mathbb{E}$.

Let $\Gamma = \Gamma_{\mathbf{A},\mathbf{B}}$ denote the random
alignment graph of the sequence-pair $( \mathbf{A},\mathbf{B} )$. In the
appropriate limit, if the alignment parameters are in the so-called
``logarithmic phase'' \cite{13,14} (i.e., if the optimal global alignment
score of long random sequences has a negative score), the random local
maximum $\tilde{M}_{m,n}$ follows an approximate Gumbel extreme value
distribution with ``scale parameter'' $\lambda$ and ``pre-factor''
$K$ \cite{15,16},
%
\begin{equation}\label{eq11}
\mathbb{P}( \tilde{M}_{m,n} > y ) \approx 1 - \exp[ - Kmn\exp( -
\lambda y ) ].
\end{equation}

\subsection{The dynamic programming algorithm for global sequence
alignment}\label{sec22}

For affine gaps $w_{g} = \Delta_{0} + \Delta_{1}g$, the global score
$S_{i,j}$ is calculated with the recursion
%
\begin{equation}\label{eq12}
S_{i,j} = \max\{ S_{i - 1,j - 1},I_{i - 1,j - 1},D_{i - 1,j - 1} \} +
s( A_{i},B_{j} ),
\end{equation}
where
\[
I_{i,j} = \max\{ S_{i,j - 1} - \Delta_{0} - \Delta_{1},I_{i,j - 1}
- \Delta_{1},D_{i,j - 1} - \Delta_{0} - \Delta_{1} \},
\]
$D_{i,j} = \max\{S_{i - 1,j} - \Delta_{0} - \Delta_{1},D_{i - 1,j} - \Delta_{1} \}$ and
boundary conditions $S_{0,0} = 0,I_{0,0} = D_{0,0} = - \infty , D_{g,0} =
I_{0,g} = - \Delta_{0} - \Delta_{1}g, S_{g,0} = S_{0,g} = I_{g,0} =
D_{0,g} = - \infty$ for $g > 0$ \cite{17}. The three array names, $S, I$,
and $D$, are mnemonics for ``substitution,'' ``insertion'' and
``deletion.'' If ``$\Delta$'' denotes a gap character, the corresponding
alignment letter-pairs $( a,b ), ( \Delta,b )$ and $( a,\Delta)$
correspond to the operations for editing sequence $\mathbf{A}$ into
sequence $\mathbf{B}$ \cite{18}.

\subsection{Previous methods for estimating the BLAST $p$-value}\label{sec23}

If $w_{g} \equiv\infty$ identically, so northward and eastward (gap)
edges are disallowed in an optimal alignment path, a rigorous proof of
(\ref{eq11}) yields analytic formulas for the Gumbel parameters
$\lambda$ and $K$ \cite{14}. For gapped local alignment, rigorous results are
sparse, although some approximate analytical studies are extant \cite{7,19,20,21}.
The prevailing approach therefore estimates $\lambda$ and $K$
from simulations \cite{22,23}. Because $\lambda$ is an exponential rate, it
dominates $K$'s contribution to the BLAST $p$-value. Most studies
therefore (including the present one) have focused on $\lambda$. (Note,
however, some recent progress on the real-time estimation of $K$ \cite{6}.)
Typically, current applications require a 1--4\% relative error
in $\lambda$; 10--20\%, in $K$ \cite{23}. The characteristics of the relevant
sequence database determine the actual accuracies required, however,
making approximations with controlled error and of arbitrary accuracy
extremely desirable in practice.

Storey and Siegmund \cite{7} approximate $\lambda$ (with neither controlled
errors nor arbitrary accuracy) as
%
\begin{equation}\label{eq13}
\tilde{\lambda} \approx\lambda^{*} - 2( \mu^{*} )^{ - 1}\Lambda e^{ -
\lambda ^{*}\Delta _{0}} / ( e^{\lambda ^{*}\Delta _{1}} - 1 ),
\end{equation}
where $\sum_{( a,b )} p_{a}p'_{b}\exp[ \lambda^{*}s( a,b ) ] = 1$ [so
$\lambda^{*}$ is the so-called ``ungapped lambda,'' for $\Delta( g )
\equiv\infty$] and $\mu^{*}: = \sum_{( a,b )} s( a,b )p_{a}p'_{b}\exp[
\lambda^{*}s( a,b ) ]$. In (\ref{eq13}), $\Lambda$ is an upper bound
for an infinite sequence of constants defined in terms of gap lengths in
a random alignment.

Many other studies have used local alignment simulations to estimate
BLAST $p$-values, for example, Chan \cite{9} used importance sampling and a
mixture distribution. Some rigorous results \cite{24} are also extant for the
so-called ``island method'' \cite{22,25}, which yields maximum likelihood
estimates of $\lambda$ and $K$ from a Poisson process associated with
local alignments exceeding a threshold score \cite{23,26}.

Large deviations arguments \cite{13,27} support the common belief that
global alignment can estimate $\lambda$ for local alignment through the
equation $\lambda =\break - \lim_{y \to \infty} y^{ - 1}\ln\mathbb{P}\{ M \ge y
\}$. For a fixed error, global alignment typically requires less
computational effort than local alignment. For example, one early study
\cite{10} used importance sampling based on trial distributions $\mathbb{Q}$
from a hidden Markov model.

The study demonstrated that the global alignment equation $\mathbb{E}[
\exp( \lambda S_{n,n} ) ] = 1$ estimated $\lambda$ with only $O( n^{ -
1} )$ error \cite{8}. (Recall that ``$\mathbb{E}$'' denotes the\break expectation
corresponding to the random letters model.) The equation\break $\mathbb{E}[
\exp( \lambda M_{m} ) ] = \mathbb{E}[ \exp( \lambda M_{n} ) ]$ $(m \ne
n)$, suggested by heuristic modeling with Markov additive processes
(MAPs) \cite{28,29}, improved the error substantially, to $O(
\varepsilon^{n} )$ \cite{5}.

The next subsection shows how the MAP heuristic can improve the
efficiency of importance sampling even further, with its renewal
structure. The next subsection gives the relevant parts of the MAP
heuristic.

\subsection{The Markov additive process heuristic}\label{sec24}

The rigorous theory of MAPs appears elsewhere \cite{28,29}. Because the MAP
heuristics given below parallel a previous publication \cite{5}, we present
only informal essentials.

Consider a finite Markov-chain state-space $\mathfrak{J}$, containing $\#\mathfrak{J}$
elements. Without loss of generality, $\mathfrak{J} = \{ 1,\ldots,\# \mathfrak{J}\}$. Until
further notice, all vectors are row vectors of dimension $\#\mathfrak{J}$; all
matrices, of dimension $( \#\mathfrak{J}) \times( \# \mathfrak{J} )$. A MAP can be defined in
terms of a time-homogenous Markov chain (MC) $( J_{n} \in \mathfrak{J}\dvtx n = 0,1,\ldots)$
and a $( \# \mathfrak{J} ) \times( \# \mathfrak{J} )$ matrix of real random variates $\|
Z_{i,j} \|$. Let the MC have transition matrix $\mathbf{P} = \| p_{i,j}
\|$, so $p_{i,j} = \mathbb{P}( J_{n} = j| J_{n - 1} = i)$. Let the
stationary distribution of the MC be $\bolds{\pi}$, assumed strictly
positive and satisfying both $\bolds{\pi} \mathbf{P} = \bolds{\pi}$
and $\bolds{\pi} \mathbf{1}^{t} = 1$, where $\mathbf{1}^{t}$ denotes the $( \#
\mathfrak{J}
) \times 1$ column vector whose elements are all $1$.

As usual, let $\mathbb{P}_{\bolds{\gamma}}$ and
$\mathbb{E}_{\bolds{\gamma}}$ be the probability measure and
expectation corresponding to an initial state $J_{0}$ with
distribution $\bolds{\gamma}$; $\mathbb{P}_{i}$ and $\mathbb{E}_{i}$, to
an initial state $J_{0} = i$; and $\mathbb{P}_{\bolds{\pi}}$
and $\mathbb{E}_{\bolds{\pi}}$, to an initial state in the equilibrium
distribution $\bolds{\pi}$.

Run the MC $( J_{n} )$, and take its succession of states as given.
Consider the following sequence $( Y_{n} \in\mathbb{R}\dvtx n = 0,1,\ldots)$ of
random variates. Define $Y_{0}: = 0$. For $n = 1,2,\ldots,$ let the $( Y_{n}
)$ be conditionally independent, with distributions determined by the
transition $J_{n - 1} \to J_{n}$ of the Markov chain as follows. If
$J_{n - 1} = i$ and $J_{n} = j$, the value of $Y_{n}$ is chosen randomly
from the distribution of~$Z_{i,j}$. (Thus, if $J_{m - 1} = J_{n - 1} = i$
and $J_{m} = J_{n} = j, Y_{m}$ and $Y_{n}$ share the distribution
of~$Z_{i,j}$, although independence permits randomness to give them
different values.)

The random variates of central interest are the sums $T_{n} = \sum_{m =
0}^{n} Y_{m}$ $(n = 0,1,\ldots)$ and the maximum $M: = \max_{n \ge 0}T_{n}$.
To exclude trivial distributions for~$M$ (i.e., $M = 0$ a.s. and $M =
\infty$ a.s.), make two assumptions: (1) $\mathbb{E}_{\bolds{\pi}} Y_{1}
< 0$; and (2)~there is some $m$ and state $i$ such that
%
\begin{equation}\label{eq14}
\mathbb{P}_{i}\bigl\{ \min\{ T_{k}\dvtx k = 1,\ldots,m \} > 0;J_{m} = i,J_{j} \ne
i\mbox{ for }j = 1,\ldots,m - 1 \bigr\} > 0.\hspace*{-28pt}
\end{equation}
Consider the sequence $( T_{n} )$, its SALEs $\beta( 0 ): = 0$ and $\beta(
k + 1 ): = \min\{ n > \beta( k )\dvtx T_{n} > T_{\beta ( k )} \}$, and its
SALE scores $T_{\beta ( k )}$. For brevity, let $\beta : = \beta( 1 )$.
Note that $M = T_{\beta ( k )}$ for some $k \in\{ 0,1,\ldots\}$. In a MAP,
$( J_{\beta ( k )},T_{\beta ( k )} )$ forms a defective Markov renewal
process.

Now, define the matrix $\mathbf{L}_{\theta} : = \| \mathbb{E}_{i}[ \exp(
\theta T_{\beta} )$; $J_{\beta} = j,\beta < \infty] \|$. The
Perron--Frobenius theorem \cite{29}, page 25, shows that $\mathbf{L}_{\theta}$
has a strictly dominant eigenvalue $\rho( \theta) > 0$ [i.e., $\rho(
\theta$) is the unique eigenvalue of greatest absolute value]. Moreover,
$\rho( \theta)$ is a convex function \cite{30}, and because $\mathbf{L}_{0}$
is substochastic, $\rho( 0 ) < 1$. The two assumptions above
(\ref{eq14}) ensure that $M: = \max_{n \ge 0}T_{n}$ has a nontrivial
distribution and that $\rho( \lambda) = 1$ for some unique $\lambda > 0$.

The notation intentionally suggests a heuristic analogy between MAPs and
global alignment. Identify the Markov chain states $J_{n}$ in the MAP
with the rectangle $[ 0,n ] \times[ 0,n ]$
of $\Gamma_{\mathbf{A},\mathbf{B}}$, and identify the sum $T_{n}$ in the
MAP with the edge maximum $M_{n}$ in global alignment. In the following,
therefore, the identification leads to $M_{n}$ replacing $T_{n}$ in the
MAP formulas. In particular, the MAP heuristic identifies the Gumbel
scale parameter in (\ref{eq11}) with the root $\lambda > 0$ of the
equation $\rho( \lambda) = 1$. Although the heuristic analogy between
MAPs and global alignment is in no way precise or rigorous, it has
produced useful results \cite{5}.

The details of why the MAP heuristic works so well are presently
obscure, although some additional motivation appears in an heuristic
calculation related to~$\lambda$ \cite{31}. The calculation takes the limit
of nested successively wider semi-infinite strips, each strip having
constant width and propagating itself northeastward in the alignment
graph $\Gamma_{\mathbf{A},\mathbf{B}}$. The successive northeast
boundaries of the propagation are states in an ergodic MC. MAPs
therefore might rigorously justify the heuristic calculation.

\section{Methods}\label{sec3}

\subsection{A novel equation for $\lambda$}\label{sec31}

From the definition of $\mathbf{L}_{\theta}$ in a MAP, if the Markov
chain $\{ J_{n} \}$ starts in a state $J_{0}$ with distribution
$\bolds{\gamma}$ (with $M_{n}$ replacing $T_{n}$ in the MAP formulas),
matrix algebra applied to the concatenation of SALEs in a MAP yields
%
\begin{equation}\label{eq21}
\mathbb{E}_{\bolds{\gamma}} \bigl[ \exp\bigl( \theta M_{\beta ( k )} \bigr);\beta( k
) < \infty\bigr] = \bolds{\gamma} ( \mathbf{L}_{\theta}
)^{k}\mathbf{1}^{t}.
\end{equation}
For a MAP, equation (\ref{eq21}) is exact; but for global alignment, it has no
literal meaning. Equation (\ref{eq21}) has some consequences for the limit $k
\to\infty$, and we speculate that the consequences hold, even for global
alignment. [Note: although the sequence $( \beta( k ) )$ is a.s. finite,
the limits $k \to\infty$ below involve no contradiction or
approximation, because they are not a.s. limits.]

Define $K_{k}( \theta): = \ln\{ \mathbb{E}_{\bolds{\gamma}} [ \exp(
\theta M_{\beta ( k )} );\beta( k ) < \infty] \}$. In (\ref{eq21}), a
spectral (eigenvalue) decomposition of the matrix $\mathbf{L}_{\theta}$
\cite{32} shows that
%
\begin{equation}\label{eq22}
K_{k}( \theta) = k\ln\{ \rho( \theta) \} + c_{0} + O( \varepsilon^{k}
),
\end{equation}
where $0 \le\varepsilon < 1$ is determined by the magnitude of the
subdominant eigenvalue of~$\mathbf{L}_{\theta}$, and $c_{0}$ is a
constant independent of $\theta$ and $k$.

For $k' - k > 0$ fixed, we can accelerate the convergence in
(\ref{eq22}) as $k \to\infty$ by differencing
%
\begin{equation}\label{eq23}
K_{k'}( \theta) - K_{k}( \theta) = ( k' - k )\ln\{ \rho( \theta) \} +
O( \varepsilon^{k} ).
\end{equation}
Let $\lambda_{k',k}$ denote the root of (\ref{eq23}) after dropping
the error term $O( \varepsilon^{k} )$. Because $\rho( \lambda) = 1$,
Taylor approximation around $\lambda$ yields $\ln\{ \rho( \lambda_{k',k}
) \} \approx\rho '( \lambda)( \lambda_{k',k} - \lambda)$, so
(\ref{eq23}) becomes
%
\begin{equation}\label{eq24}
( k' - k )\rho '( \lambda)( \lambda_{k',k} - \lambda) = O(
\varepsilon^{k} ),
\end{equation}
that is, with $k' - k$ fixed, $\lambda_{k',k}$ converges geometrically to
$\lambda$ as the SALE index $k \to\infty$.

The initial state $\bolds{\gamma}$ of global alignment has a
deterministic distribution, namely the origin $( 0,0 )$. Equation (\ref{eq23})
for $\theta = \lambda$ therefore becomes
%
\begin{equation}\label{eq25}
\mathbb{E}\bigl[ \exp\bigl( \lambda M_{\beta ( k' )} \bigr);\beta( k' ) < \infty\bigr] =
\mathbb{E}\bigl[ \exp\bigl( \lambda M_{\beta ( k )} \bigr);\beta( k ) < \infty\bigr]
\end{equation}
after dropping the geometric error $O( \varepsilon^{k} )$. Let
$\hat{\lambda} _{k',k}$ be the root of~(\ref{eq25}).

\subsection{The trial distribution for importance sampling}\label{sec32}

In (\ref{eq25}), crude Monte Carlo simulation generating random
sequence-pairs with the identical letters model $\mathbb{P}$ is
inefficient for the following reason.
When practical alignment scoring systems are used,
$\mathbb{P}\{ \beta( k ) < \infty\} < 1$ for $k \ge
1$. For, example, the BLAST defaults (scoring matrix BLOSUM62, gap
penalty $w_{g} = 11 + g$, and Robinson--Robinson letter frequencies),
$\mathbb{P}\{ \beta( 4 ) < \infty\} \approx 0.047$, so only about 1 in
20 crude Monte Carlo simulations generate a fourth ladder point.
Empirically in our importance sampling, however, Gumbel parameter
estimation seemed most efficient when the stopping time corresponded to
$\beta( 4 )$ (see below).

Importance sampling requires a trial distribution to determine
$\hat{\lambda} _{k',k}$ from (\ref{eq25}). By editing one sequence
into another, a Markov chain model borrowed directly from a previous
study \cite{10} generates random sequence alignments, as follows.

Consider a Markov state space consisting of the set of alignment
letter-pairs~$\bar{\mathfrak{L}}^{2}$, where $\bar{\mathfrak{L}}: = \mathfrak{L} \cup\{ \Delta\}$,
``$\Delta$'' being a character representing gaps. The ordered pair $(
\Delta,\Delta)$ has probability 0, so a succession of Markov states
corresponds to a global sequence alignment (see Figure \ref{fig1}), that is, to a
path in the alignment graph $\Gamma_{\mathbf{A},\mathbf{B}}$. Ordered
pairs other than $( \Delta,\Delta)$ fall into three sets, corresponding
to edit operations following (\ref{eq12}): $S: = \mathfrak{L} \times \mathfrak{L}$
[substitution, a bioinformatics term implicitly including identical
letter-pairs $( a,a )$], $I: = \{ \Delta\} \times \mathfrak{L}$ (insertion); and $D: =
\mathfrak{L} \times\{ \Delta\}$ (deletion). The sets $S, I$ and $D$ form ``atoms''
of the MC \cite{33}, page 203, as follows. (By definition, each atom of a MC
is a set of all states with identical outgoing transition probabilities.)

From the set $S$, the transition probability to $( a,b )$
is $t_{S,S}q_{a,b}$; to $( \Delta,b ),t_{S,I}p'_{b}$; and to $(
a,\Delta),t_{S,D}p_{a}$. From the set $I$, the transition probability to
$( a,b )$ is $t_{I,S}q_{a,b}$; to $( \Delta,b ),t_{I,I}p'_{b}$; and to $(
a,\Delta),t_{I,D}p_{a}$. From the set $D$, the transition probability to
$( a,b )$ is $t_{D,S}q_{a,b}$; to $( \Delta,b ),t_{D,I}p'_{b}$; and to $(
a,\Delta),t_{D,D}p_{a}$. Transition probabilities sum to 1, so the
following restrictions apply: $\sum_{a,b \in \mathfrak{L}} q_{a,b} = 1,
\sum_{b \in \mathfrak{L}} p'_{b} = 1, \sum_{a \in \mathfrak{L}} p_{a} =
1, t_{S,S} + t_{S,I} + t_{S,D} = 1$ (transit from the substitution
atom), $t_{D,D} + t_{D,S} + t_{D,I} = 1$ (transit from the deletion
atom) and $t_{I,I} + t_{I,S} + t_{I,D} = 1$ (transit from the insertion
atom). Usually in practice, the term $t_{I,D} = 0$, to disallow
insertions following a deletion. Our formulas retain the term, to
exploit the resulting symmetry later.

In the terminology of hidden Markov models, $S, I, D$ are hidden Markov
states. $t_{i,j}$ for $i,j \in\{ S,I,D \}$ are transition probabilities
and $q_{a,b}, p'_{b}, p_{a}$ for $a,b \in\mathfrak{L}$ are emission
probabilities from the state $S, I, D$, respectively.

As described elsewhere \cite{10}, numerical values for the Markov
probabilities can be determined from the scores $s( a,b )$ and the gap
penalty $w_{g}$. Note that the values are selected for statistical
efficiency, although many other values also yield unbiased estimates for
$\lambda$ in the appropriate limit.

\subsection{Importance sampling weights and stopping times}\label{sec33}

To establish notation, and to make connections to the \hyperref[app]{Appendix} and its
mapping theorem, note that the MC above can be supported on a
probability space $( \Omega,\mesf{F},\mathbb{Q} )$, where each $\omega =
( \pi,\mathbf{A},\mathbf{B} ) \in\Omega$ is an ordered triple. Here,
$\pi$ is an infinite path starting at the origin in the alignment
graph $\Gamma_{\mathbf{A},\mathbf{B}}; \mesf{F}$ is the set generated by
cylinder sets in $\Omega$ (here, cylinder sets essentially consist of
some finite path and the corresponding pair of subsequences); and
$\mathbb{Q}$ is the MC probability distribution described above, started
at the atom $S$, with expectation operator $\mathbb{E}_{\mathbb{Q}}$.

Let $N$ be any stopping time for the sequence $( M_{n}\dvtx n = 0,1,\ldots)$ of
edge maxima for $\Gamma_{\mathbf{A},\mathbf{B}}$ (i.e., the sequence $\{
M_{0},\ldots,M_{n} \}$ determines whether $N \le n$ or not). Because
$M_{n}$ is determined by $( \mathbf{A}[ 1,n ],\mathbf{B}[ 1,n ] ), N$ is
also a stopping time for the sequence $\{ ( \mathbf{A}[ 1,n ],\mathbf{B}[
1,n ] )\dvtx n = 0,1,\ldots\}$. The stopping time of main interest here is $N =
\beta( k )$, the $k$th ladder index of $( M_{n} )$, where $k \ge 1$ is
arbitrary. (As further motivation for the mapping theorem in the
\hyperref[app]{Appendix}, other stopping times of possible interest include, for
example, $N = n$, a fixed epoch \cite{8}, and $N = \beta( K_{y} )$, where
$\beta( K_{y} ) = \inf\{ n\dvtx M_{n} \ge y \}$ is the index of first
ladder-score outside the interval $( 0,y ).)$

To use the mapping theorem, introduce the probability space $( \Omega
'',\mesf{F}'',\mathbb{P} )$, where each $\omega '' = (
\mathbf{A},\mathbf{B} ) \in\Omega ''$ is an ordered pair. Here,
$\mathbf{A}$ and $\mathbf{B}$ are sequences, $\mesf{F}''$ is the set
generated by all cylinder sets in $\Omega ''$ (i.e., sets corresponding
to pairs of finite subsequences) and $\mathbb{P}( A'' ) = \prod_{k =
1}^{i} p_{A_{k}} \prod_{k = 1}^{j} p'_{B_{k}}$, if the cylinder set
$A''$ corresponds to the subsequence pair $( \mathbf{A}[ 1,i
],\mathbf{B}[ 1,j ] )$. Given $N$, the theory of stopping times \cite{29}, page
414, can be used to construct a discrete probability space $( \Omega
',\mesf{F}',\mathbb{P} )$, where each event $\omega ' \in\Omega '$ is a
finite-sequence pair $\omega ' = ( \mathbf{A}[ 1,N ],\mathbf{B}[ 1,N ] ),
\mesf{F}'$ is the set of all subsets of $\Omega '$ and $\mathbb{P}(
\omega ' ) = \prod_{k = 1}^{N( \omega ' )} p_{A_{k}} \prod_{k = 1}^{N(
\omega ' )} p'_{B_{k}}$.

\begin{figure}[b]

\includegraphics{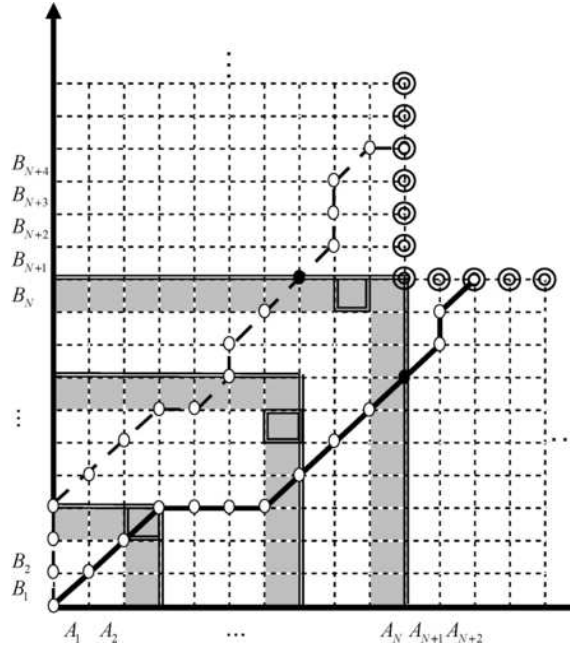}

\caption{Two examples of alignment path $\pi$ generated by a Markov
chain. As in Figure \protect\ref{fig1}, the shading and the double lines indicate squares
where a vertex (surrounded by double lines) generated an SALE. The SALEs
determine the stopping time $N = \beta( 3 )$. In Figure \protect\ref{fig2}, the first SALE
is determined by the score at the vertex $( 3,3 );$ the second SALE, the
vertex $( 7,6 );$ the third SALE, the vertex $( 9,10 )$. Therefore, $N =
\beta( 3 ) = 10$. The vertical ray $\mesf{I}_{N,N}$ and the horizontal
ray  $\mesf{D}_{N,N}$ are indicated by double circles. The lower path
$\pi$ (solid line) ends at $( N + 2,N )$ with a final transition to $S$;
the upper path $\pi$ (long-dashed line), at $( N,N + 4 )$ with a final
transition to $D$. The closed vertices indicate intersection with the
square corresponding to $\omega ' = ( \omega '_{\mathbf{A}},\omega
'_{\mathbf{B}} ) = ( \mathbf{A}[ 1,N ],\mathbf{B}[ 1,N ] )$.}
\label{fig2}
\end{figure}

Let $\mesf{I}_{m,n}: = \{ ( i,j )\dvtx i = m,j \ge n \}$
and $\mesf{D}_{m,n}: = \{ ( i,j )\dvtx i \ge m,j = n \}$. Define the
function $f\dvtx\omega\mapsto\omega '$, where $\omega = (
\pi,\mathbf{A},\mathbf{B} )$ and $\omega ' = ( \omega
'_{\mathbf{A}},\omega '_{\mathbf{B}} ): = ( \mathbf{A}[ 1,N
],\break\mathbf{B}[ 1,N ] )$. Then, $\omega\in f^{ - 1}( \omega ' )$, if and
only if the path $\pi$ hits the set $\mesf{I}_{N,N}
\cup\mesf{D}_{N,N}$ at $( i,j )$, so that $\mathbf{A}[ 1,N ] = \omega
'_{\mathbf{A}}$ and $\mathbf{B}[ 1,N ] = \omega '_{\mathbf{B}}$ (see
Figure \ref{fig2}).

Empirically, our simulations satisfied $\mathbb{Q}\{ \beta( k ) <
\infty\} = 1$, and we speculate that our application therefore satisfies
the hypothesis $\mathbb{Q}H = 1$ of the \hyperref[app]{Appendix}. According to the
\hyperref[app]{Appendix}, the reciprocal importance sampling weight $1 / W( \omega) =
\sum_{\omega _{0} \in f^{ - 1}\{ f( \omega ) \}} \mathbb{Q}( \omega_{0}
) /\mathbb{P}f( \omega)$ depends on the sum over all possible Markov
chain realizations $\omega_{0} \in f^{ - 1}( \omega ' )$. Dynamic
programming computes the sum efficiently, as follows.

Let the ``transition'' $T$ represent any element of $\{ S,I,D \}$
[substitution $( a_{i},b_{j} )$, insertion $( \Delta,b_{j} )$, or
deletion $( a_{i},\Delta)$]. Fix any particular pair $(
\mathbf{A},\mathbf{B} )$ of infinite sequences, which fixes $N = \beta( k
)$. To set up a recursion for dynamic programming, consider the
following set of events $\mesf{E}_{i,j}^{T}$, defined for $T \in\{ S,I,D
\}$ and \mbox{$\min\{ i,j \} \le N$}, and illustrated in Figure \ref{fig2}.\vspace*{2pt} Let
$\mesf{E}_{i,j}^{T}$ be the event consisting of all $\omega$ yielding a
path $\pi$ whose final transition is $T$ and which corresponds to the
subsequences: (1) $\mathbf{A}[ 1,i ]$ and $\mathbf{B}[ 1,j ]$ for $0 \le
i,j \le N$; (2) $\mathbf{A}[ 1,i ]$ and $\mathbf{B}[ 1,N ]$ for $0 \le N = j
\le i$; and (3) $\mathbf{A}[ 1,N ]$ and $\mathbf{B}[ 1,j ]$ for $0 \le N =
i \le j$. Define $Q_{i,j}^{T}: = \mathbb{Q}( \mesf{E}_{i,j}^{T} )$
and $Q_{i,j}: = Q_{i,j}^{S} + Q_{i,j}^{I} + Q_{i,j}^{D}$. (Note: in the
following, $T \in\{ S,I,D \}$ is always a superscript, never an
exponent.)

For brevity, let $\tilde{q}_{i,j} = q_{A_{i},B_{j}}$ for $0 \le i,j \le
N; \tilde{q}_{i,j} = \sum_{( a \in \mathfrak{L} )} q_{a,B_{j}}$ for $0
\le j \le N < i; \tilde{q}_{i,j} = \sum_{( b \in \mathfrak{L} )}
q_{A_{i},b}$ for $0 \le i \le N < j$; and $\tilde{q}_{i,j} = 1$
otherwise. Let $\tilde{p}'_{j} = p'_{B_{j}}$ for $0 \le j \le N$; and $1$
otherwise. Finally, Let $\tilde{p}_{i} = p_{A_{i}}$ for $0 \le i \le N$;
and~$1$ otherwise. Because every path into the vertex $( i,j )$ comes
from one of three vertices, each corresponding to a different
transition $T \in\{ S,I,D \}$,
%
\begin{eqnarray}\label{eq26}
Q_{i,j}^{S} &=& \tilde{q}_{i,j}( t_{S,S}Q_{i - 1,j -
1}^{S} + t_{I,S}Q_{i - 1,j - 1}^{I} + t_{D,S}Q_{i - 1,j - 1}^{D} ), \nonumber\\
Q_{i,j}^{I} &=& \tilde{p}'_{j}( t_{S,I}Q_{i,j - 1}^{S} + t_{I,I}Q_{i,j -
1}^{I} + t_{D,I}Q_{i,j - 1}^{D} ), \\
Q_{i,j}^{D} &=& \tilde{p}_{i}(t_{S,D}Q_{i - 1,j}^{S} + t_{I,D}Q_{i - 1,j}^{I} + t_{D,D}Q_{i - 1,j}^{D}
)\nonumber
\end{eqnarray}
with boundary conditions $Q_{0,0}^{S} = 1, Q_{0,0}^{I} = Q_{0,0}^{D} = 0,
Q_{g,0}^{S} = Q_{0,g}^{S} = Q_{g,0}^{I} = Q_{0,g}^{D} = 0, Q_{0,g}^{I} =
p'_{B_{1}}\cdots p'_{B_{g}}t_{S,I}( t_{I,I} )^{g - 1}$ and $Q_{g,0}^{D} =
p_{A_{1}}\cdots p_{A_{g}}t_{S,D}\times\break( t_{D,D} )^{g - 1}$ $(g > 0)$.

Recall that $\omega = ( \pi,\mathbf{A},\mathbf{B} ) \in f^{ - 1}( \omega
' )$, if and only if the path $\pi$ hits the set $\mesf{I}_{N,N}
\cup\mesf{D}_{N,N}$ at $( i,j )$, so that $\mathbf{A}[ 1,N ] = \omega
'_{\mathbf{A}}$ and $\mathbf{B}[ 1,N ] = \omega '_{\mathbf{B}}$. Thus,
%
\begin{equation}\label{eq27}\qquad
\sum_{\omega \in f^{ - 1}( \omega ' )} \mathbb{Q}( \omega) = -
Q_{N,N}^{S} + \sum_{j = N}^{\infty} ( Q_{N,j}^{S} + Q_{N,j}^{D} ) +
\sum_{i = N}^{\infty} ( Q_{i,N}^{S} + Q_{i,N}^{I} ).
\end{equation}
To turn (\ref{eq26}) into a recursion for importance sampling weights,
define $P_{i}: = p_{A_{1}}\cdots p_{A_{\min \{ i,N \}}} =
\tilde{p}_{1}\cdots\tilde{p}_{i}$ and $P'_{j}: = p'_{B_{1}}\cdots p'_{B_{\min \{
j,N \}}} = \tilde{p}'_{1}\cdots\tilde{p}'_{j}$, and let $W_{i,j}^{T}: =
Q_{i,j}^{T}/( P_{i}P'_{j} )$\vspace*{2pt} $(T \in\{ S,I,D \})$. Let $r_{i,j} =
\tilde{q}_{i,j}/( \tilde{p}_{i}\tilde{p}'_{j} )$. For future reference,
define $r_{ \bullet ,j}: = r_{i,j}$ for $0 \le j \le N < i$ and $r_{i,
\bullet} : = r_{i,j}$ for $0 \le i \le N < j$. Note that $r_{ \bullet
,j}$ is independent of $i$, and $r_{i, \bullet}$ is independent of $j$.
Equation (\ref{eq26}) yields
%
\begin{eqnarray}\label{eq28}
W_{i,j}^{S} &=& r_{i,j}( t_{S,S}W_{i - 1,j - 1}^{S} +
t_{I,S}W_{i - 1,j - 1}^{I} + t_{D,S}W_{i - 1,j - 1}^{D} ), \nonumber\\
W_{i,j}^{I} &=& t_{S,I}W_{i,j - 1}^{S} + t_{I,I}W_{i,j - 1}^{I} + t_{D,I}W_{i,j -
1}^{D}, \\
W_{i,j}^{D} &=& t_{S,D}W_{i - 1,j}^{S} + t_{I,D}W_{i - 1,j}^{I}
+ t_{D,D}W_{i - 1,j}^{D}\nonumber
\end{eqnarray}
with boundary conditions $W_{0,0}^{S} = 1, W_{0,0}^{I} = W_{0,0}^{D} = 0,
W_{g,0}^{S} = W_{0,g}^{S} = W_{0,g}^{I} = W_{g,0}^{D} = 0, W_{0,g}^{I} =
t_{S,I}( t_{I,I} )^{g - 1}$ and $W_{g,0}^{D} = t_{S,D}( t_{D,D} )^{g -
1}$ $(g > 0)$. Because\vspace*{1pt} of (\ref{eq27}), the importance sampling weight
$W: = W( \omega)$ satisfies
%
\begin{eqnarray}\label{eq29}
\frac{1}{W} &=& \frac{\sum_{\omega _{0} \in f^{ - 1}\{ f( \omega ) \}}
\mathbb{Q}( \omega _{0} )} {\mathbb{P}f( \omega )}\nonumber\\[-8pt]\\[-8pt]
&=& - W_{N,N}^{S} +
\sum_{j = N}^{\infty} ( W_{N,j}^{S} + W_{N,j}^{D} ) + \sum_{i =
N}^{\infty} ( W_{i,N}^{S} + W_{i,N}^{I} ).\nonumber
\end{eqnarray}

Because $r_{i,j} = r_{ \bullet ,j}$ $(0 \le j \le N < i)$ and $r_{i,j} =
r_{i, \bullet}$ $(0 \le i \le N < j)$, only a finite number of recursions
are needed to compute the infinite sums in (\ref{eq29}), as follows.
For $T \in\{ S,I,D \}$, define $\tilde{U}_{i}^{T}: = U_{i,N}^{T}$,
where $U_{m,n}^{T}: = \sum_{j = n}^{\infty} W_{m,j}^{T}$. Likewise,
define $\tilde{V}_{j}^{T}: = V_{N,j}^{T}$, where $V_{m,n}^{T}: = \sum_{i =
m}^{\infty} W_{i,n}^{T}$. Equation (\ref{eq29}) becomes
%
\begin{equation}\label{eq210}
\frac{1}{W} = - W_{N,N}^{S} + \tilde{U}_{N}^{S} + \tilde{U}_{N}^{D} +
\tilde{V}_{N}^{S} + \tilde{V}_{N}^{I}.
\end{equation}

Note that $U_{i,j - 1}^{T} - U_{i,j}^{T} = W_{i,j - 1}^{T}$. To
determine $\tilde{U}_{N}^{T}$, summation of (\ref{eq28}) for $0 \le i
\le N < j$ yields
%
\begin{eqnarray}\label{eq211}\qquad
U_{i,j}^{S} &=& r_{i, \bullet} ( t_{S,S}U_{i - 1,j -
1}^{S} + t_{I,S}U_{i - 1,j - 1}^{I} + t_{D,S}U_{i - 1,j - 1}^{D} )\nonumber\\
& =&
U_{i,j - 1}^{S} - W_{i,j - 1}^{S}, \nonumber\\[-8pt]\\[-8pt]
U_{i,j}^{I} &=& t_{S,I}U_{i,j -
1}^{S} + t_{I,I}U_{i,j - 1}^{I} + t_{D,I}U_{i,j - 1}^{D} = U_{i,j -
1}^{I} - W_{i,j - 1}^{I} ,\nonumber\\
U_{i,j}^{D} &=& t_{S,D}U_{i - 1,j}^{S} +
t_{I,D}U_{i - 1,j}^{I} + t_{D,D}U_{i - 1,j}^{D}.\nonumber
\end{eqnarray}
Elimination of $U_{i,j}^{T}$ for $j = N + 1$ and $i = 1,\ldots,N$ in the
first two equations yields
%
\begin{eqnarray}\label{eq212}
U_{i,N}^{S} &=& r_{i, \bullet} ( t_{S,S}U_{i - 1,N}^{S} +
t_{I,S}U_{i - 1,N}^{I} + t_{D,S}U_{i - 1,N}^{D} ) + W_{i,N}^{S},
\nonumber\\
U_{i,N}^{I} &=& t_{S,I}U_{i,N}^{S} + t_{I,I}U_{i,N}^{I} +
t_{D,I}U_{i,N}^{D} + W_{i,N}^{I}, \\
U_{i,N}^{D} &=& t_{S,D}U_{i -
1,N}^{S} + t_{I,D}U_{i - 1,N}^{I} + t_{D,D}U_{i - 1,N}^{D},\nonumber
\end{eqnarray}
that is,
%
\begin{eqnarray}\label{eq213}
\tilde{U}_{i}^{S} &=& r_{i, \bullet} ( t_{S,S}\tilde{U}_{i
- 1}^{S} + t_{I,S}\tilde{U}_{i - 1}^{I} + t_{D,S}\tilde{U}_{i - 1}^{D} )
+ W_{i,N}^{S}, \nonumber\\
\tilde{U}_{i}^{I} &=& ( 1 - t_{I,I} )^{ - 1}(
t_{S,I}\tilde{U}_{i}^{S} + t_{D,I}\tilde{U}_{i}^{D} + W_{i,N}^{I} ), \\
\tilde{U}_{i}^{D} &=& t_{S,D}\tilde{U}_{i - 1}^{S} + t_{I,D}\tilde{U}_{i -
1}^{I} + t_{D,D}\tilde{U}_{i - 1}^{D}\nonumber
\end{eqnarray}
with initial values $\tilde{U}_{0}^{S} = \tilde{U}_{0}^{D} = 0$
and $\tilde{U}_{0}^{I} = ( 1 - t_{I,I} )^{ - 1}W_{0,N}^{I} = ( 1 -
t_{I,I} )^{ - 1}\times\break t_{S,I}( t_{I,I} )^{N - 1}$. Compute (\ref{eq213})
recursively for $i = 1,\ldots,N$.

Similarly, reflect through $i = j$ to derive
%
\begin{eqnarray}\label{eq214}
\tilde{V}_{j}^{S} &=& r_{ \bullet ,j}( t_{S,S}\tilde{V}_{j
- 1}^{S} + t_{D,S}\tilde{V}_{j - 1}^{D} + t_{I,S}\tilde{V}_{j - 1}^{I} )
+ W_{N,j}^{S}, \nonumber\\
\tilde{V}_{j}^{I} &=& t_{S,I}\tilde{V}_{j - 1}^{S} +
t_{D,I}\tilde{V}_{j - 1}^{D} + t_{I,I}\tilde{V}_{j - 1}^{I}, \\
\tilde{V}_{j}^{D} &=& ( 1 - t_{D,D} )^{ - 1}( t_{S,D}\tilde{V}_{j}^{S} +
t_{I,D}\tilde{V}_{j}^{I} + W_{N,j}^{D} ) \nonumber
\end{eqnarray}
with initial values $\tilde{V}_{0}^{S} = \tilde{V}_{0}^{I} = 0$
and $\tilde{V}_{0}^{D} = ( 1 - t_{D,D} )^{ - 1}W_{N,0}^{D} = ( 1 -
t_{D,D} )^{ - 1}\times\break t_{S,D}( t_{D,D} )^{N - 1}$. Iterate (\ref{eq214})
for $j = 1,\ldots,N$. Substitute the results for $\tilde{U}_{N}^{S},
\tilde{U}_{N}^{D}, \tilde{V}_{N}^{S}$, and $\tilde{V}_{N}^{I}$ into
(\ref{eq210}) to compute $W$.

\subsection{Error estimates for $\hat{\lambda}_{k',k}$}\label{sec34}

Denote the indicator of an event $A$ by $\mathbb{I}A$, that is,
$\mathbb{I}A = 1$ if $A$ occurs and 0 otherwise. For a realization
$\omega$ in the simulation, define
%
%
\begin{eqnarray}\label{eq215}\qquad
h_{k,k'}( \theta) &: =& h_{k,k'}( \theta ;\omega)\nonumber\\[-8pt]\\[-8pt]
&: =& \exp\bigl( \theta M_{\beta
( k' )} \bigr)\mathbb{I}[ \beta( k' ) < \infty] - \exp\bigl( \theta M_{\beta ( k
)} \bigr)\mathbb{I}[ \beta( k ) < \infty]\nonumber
\end{eqnarray}
and let $h'_{k,k'}$ be its derivative with respect to $\theta$.

Given samples $\omega_{i}$ $(i = 1,\ldots,r)$ from the trial
distribution $\mathbb{Q}$, let $W = W( \omega_{i} )$ denote the
corresponding importance sampling weights. Because $\hat{\lambda}
_{k',k}$ is the \mbox{M-estimator} \cite{34} of the root $\lambda_{k',k}$
of $\mathbb{E}h_{k,k'}( \lambda_{k',k} ) = 0$, as $r \to\infty, \sqrt{r} (
\hat{\lambda} _{k',k} - \lambda_{k',k} )$ converges in distribution to
the normal distribution with mean $0$ and variance~\cite{34}\looseness=1
%
\begin{equation}\label{eq216}
\frac{\mathbb{E}_{\mathbb{Q}}[ h( \lambda _{k',k} )W ]^{2}}{\{
\mathbb{E}_{\mathbb{Q}}[ h'( \lambda _{k',k} )W ] \}^{2}}
\approx\frac{r^{ - 1}\sum_{1}^{r} [ h( \omega _{i};\hat{\lambda} _{k',k}
)W( \omega _{i} ) ]^{2}} {\{ r^{ - 1}\sum_{1}^{r} [ h'( \omega
_{i};\hat{\lambda} _{k',k} )W( \omega _{i} ) ] \}^{2}}.
\end{equation}

\section{Numerical study for Gumbel scale parameter}\label{sec4}

Table \ref{tab1} gives our ``best estimate'' $\bar{\lambda}$ of the Gumbel scale
parameter $\lambda$ from (\ref{eq25}) for each of the 5 options BLASTP
gives users for the alignment scoring scheme. For every scheme,
estimates $\hat{\lambda}$ derived from the first to fourth
SALEs indicated that $\hat{\lambda}$ generally is biased above the true
value $\lambda$, but that $\hat{\lambda}$ converged adequately by the
fourth SALE. The best estimate $\bar{\lambda}$ (shown in Table \ref{tab1}) is
the average of 200 independent estimates~$\hat{\lambda}$, each computed
within 1 sec from sequence-pairs simulated up to their fourth SALE.
For BLOSUM 62 and gap penalty $w_{g} = 11 + g$, the average computation
produced $1441$ sequence-pairs up to their fourth SALE within 1
second. (For results relevant to the other publicly available scoring
schemes, see Table \ref{tab1}.) The best estimates $\bar{\lambda}$ derived from
(\ref{eq25}) were within the error of the BLASTP values for $\lambda$.

\begin{table}
\caption{Best estimates $\bar{\lambda}$ for the $5$ BLASTP
alignment scoring schemes. For each scheme, we generated $200$
estimates $\hat{\lambda}$, each within a one-second computation time.
The third column gives present estimates of $\lambda$ used on the BLAST
web page (Stephen Altschul: personal communication). The BLAST values
are accurate to approximately $\pm$1\%. The fourth column gives the mean
$\bar{\lambda}$ of our $200$ estimates $\hat{\lambda}$; the fifth, the
standard error of $\bar{\lambda}$, which can be multiplied by
$\sqrt{200} \approx 14$ to give the standard error in
each $\hat{\lambda}$. The sixth column gives the average number of
sequence-pairs used to estimate each $\hat{\lambda}$. The total number
of sequence-pairs used for $\bar{\lambda}$ is $200$ times average number
of sequence-pairs. The last column gives the average sequence length
required for the fourth SALE used to estimate each $\hat{\lambda}$}
\label{tab1}
\begin{tabular*}{\tablewidth}{@{\extracolsep{\fill}}lccccd{4.0}d{2.2}@{}}
\hline
& \textbf{Gap} & \textbf{BLAST} & \textbf{Best} & \textbf{Standard} &  &
\multicolumn{1}{c@{}}{\textbf{Average}}\\
\textbf{Scoring} & \textbf{penalty} & \textbf{value} & \textbf{estimate} & \textbf{error of}
& \multicolumn{1}{c}{\textbf{Average number of}} & \multicolumn{1}{c@{}}{\textbf{sequence}}\\
\textbf{matrix} & $\bolds{w_{g}}$ &  & $\bolds{\bar{\lambda}}$
& $\bolds{\bar{\lambda}}$ & \multicolumn{1}{c}{\textbf{sequence-pairs}} & \multicolumn{1}{c@{}}{\textbf{length}}\\
\hline
BLOSUM80 & $10 + g$ & 0.299 & 0.2998 & 0.0001 & 2865 & 15.85\\
BLOSUM62 & $11 + g$ & 0.267 & 0.2679 & 0.0002 & 1441 & 27.78\\
BLOSUM45 & $14 + 2g$ & 0.195 & 0.1962 & 0.0003 & 789 & 39.23\\
PAM30 & $9 + g$ & 0.294 & 0.2956 & 0.0001 & 3593 & 9.20\\
PAM70 & $10 + g$ & 0.291 & 0.2922 & 0.0001 & 3397 & 11.49\\
\hline
\end{tabular*}
\end{table}

\begin{figure}[b]

\includegraphics{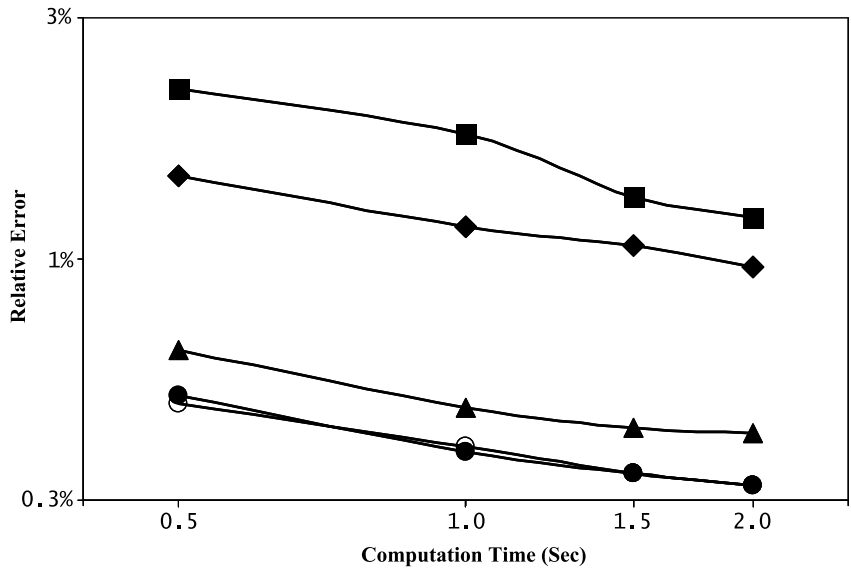}

\caption{Plot of relative errors against computation time (sec). Both
axes are in logarithmic scale. Computation time was measured on a 2.99
GHz Pentium\tsup{\textregistered}\ D CPU. Relative errors for BLOSUM45 with
$\Delta( g ) = 14 + 2g$ are shown by $\blacksquare$; BLOSUM62
with $\Delta( g ) = 11 + g$, by $\blacklozenge$; BLOSUM80 with $\Delta( g
) = 10 + g$, by \TUB; PAM70 with $\Delta( g ) = 10 + g$, by \CIB;
PAM30 with $\Delta( g ) = 9 + g$, by \CIW.}
\label{fig3}
\end{figure}

Despite having the variance formula in (\ref{eq216}) in hand, we
elected to estimate the standard error $\hat{s}_{\lambda}$ directly from
the 200 independent estimates $\hat{\lambda}$. Figure \ref{fig3} plots the
relative error $\hat{s}_{\lambda} /\bar{\lambda}$ in each individual
$\hat{\lambda}$ against the computation time, where~$\hat{s}_{\lambda}$
is the standard error of $\hat{\lambda}$. It shows that for all 5 BLASTP
online options, (\ref{eq25}) easily computed $\hat{\lambda}$ to 1--4\%
accuracy within about 0.5 seconds.

\section{Discussion}\label{sec5}

This article indicates that the scale parameter $\lambda$ of the Gumbel
distribution for local alignment of random sequences satisfies
(\ref{eq25}), an equation involving the strict ascending ladder-points (SALEs)
from global alignment, at least approximately. For standard protein
scoring systems, in fact, simulation error could account for most (if
not all) of the observed differences between values of~$\lambda$
calculated from (\ref{eq25}) and values calculated from extensive
crude Monte Carlo simulations. (The values of $\lambda$ from crude
simulation have a standard error of about $\pm$1\%.) In SALE
simulations, (\ref{eq25}) estimated $\lambda$ to 1--4\% accuracy within
0.5 second, as required by BLAST database searches over the Web. The
present study did not tune simulations much; it relied instead on
methods specific to sequence alignment to improve estimation. Many
general strategies for sequential importance sampling therefore remain
available to speed simulation. Preliminary investigations estimating the
other Gumbel parameter (the pre-factor $K$) with SALEs are encouraging,
so online estimation of the entire Gumbel distribution for arbitrary
scoring schemes appears imminent, and preliminary computer code is
already in place.

\begin{appendix}\label{app}
\section*{Appendix: A general mapping theorem for importance~sampling}

The following theorem describes an unusual type of Rao-Blackwellization
\cite{35}. Consider two probability spaces $( \Omega,\mesf{F},\mathbb{Q} )$
and $( \Omega ',\mesf{F}',\mathbb{P} )$, and a
$\mesf{F}/\mesf{F}'$-measurable function $f\dvtx\Omega\mapsto\Omega '$
(i.e., $f^{ - 1}F' \in\mesf{F}$ for every $F' \in\mesf{F}')$. Note:
$f$ is explicitly permitted to be many-to-one. Let $\mathbb{P} < <
\mathbb{Q}f^{ - 1}$ on some set $H'$ (i.e., $\mathbb{Q}f^{ - 1}G' =
0\Rightarrow \mathbb{P}G' = 0$ for any set $G' \subseteq H')$, so the
Radon--Nikodym derivative in the second line of (\ref{eq31}) below
exists. Let $H: = f^{ - 1}H'$, so for every random variate $X'$ on $(
\Omega ',\mesf{F}' )$,
%
%
\begin{eqnarray}\label{eq31}\qquad
\mathbb{E}[ X';H' ]&: =& \int_{\omega ' \in H'} X'( \omega
' )\,d\mathbb{P}( \omega ' ) \nonumber\\
&=& \int_{\omega ' \in H'} X'( \omega '
)\,d\mathbb{P}( \omega ' )\int_{\omega \in f^{ - 1}( \omega ' )}
\frac{d\mathbb{Q}( \omega )}{\int_{\omega _{0} \in f^{ - 1}( \omega ' )}
\,d\mathbb{Q}( \omega _{0} )} \nonumber\\[-8pt]\\[-8pt]
&=& \int_{\omega ' \in H'} \int_{\omega
\in f^{ - 1}( \omega ' )} X'f( \omega )\frac{d\mathbb{P}f( \omega
)}{\int_{\omega _{0} \in f^{ - 1}\{ f( \omega ) \}} \,d\mathbb{Q}( \omega
_{0} )} \,d\mathbb{Q}( \omega ) \nonumber\\
&=& \int_{\omega \in H} X'f( \omega
)\frac{d\mathbb{P}f( \omega )}{\int_{\omega _{0} \in f^{ - 1}\{ f(
\omega ) \}} \,d\mathbb{Q}( \omega _{0} )} \,d\mathbb{Q}( \omega
).\nonumber
\end{eqnarray}
Consider the application of (\ref{eq31}) to importance sampling with
target distribution~$\mathbb{P}$ and trial distribution $\mathbb{Q}$.
Assume $\mathbb{Q}H = 1$, so $H$ supports $\mathbb{Q}$. In our application
to global alignment, $H = [ \beta( k ) < \infty]
\subset\Omega$ (``$\subset$'' being strict inclusion), but we
speculate $\mathbb{Q}H = 1$.

In Monte Carlo applications, a discrete sample space $H$ is usually
available. Accordingly, the following theorem replaces the integrals in
(\ref{eq31}) by sums.

\subsection*{The mapping theorem for importance sampling} \textit{Let}
%
\begin{equation}\label{eq32}
\frac{1}{W( \omega )}: = \frac{\sum_{\omega _{0} \in f^{ - 1}\{ f(
\omega ) \}} \mathbb{Q}( \omega _{0} )} {\mathbb{P}f( \omega )}.
\end{equation}
\textit{Under the above conditions}, $r^{ - 1}\sum_{i = 1}^{r} [ X'f(
\omega_{i} )W( \omega_{i} ) ] \to\mathbb{E}[ X';H' ]$ \textit{with
probability $1$ and in mean} (\textit{with respect to} $\mathbb{Q}$), \textit{as the number
of realizations} $r \to\infty$.

The mapping theorem is an easy application of the law of large numbers
to~(\ref{eq31}).
\end{appendix}

\section*{Acknowledgments}
The authors Y. Park and S. Sheetlin contributed equally to the article.
All authors would like to acknowledge helpful discussion with Dr.
Nak-Kyeong Kim. This research was supported by the Intramural Research
Program of the NIH, National Library of Medicine.

\printaddresses

\end{document}